\documentclass[12pt]{article}
\usepackage{amsmath}
\usepackage{amsfonts}
\usepackage{epsfig}

\newcommand{\lbl}{\label}

\newcommand{\eq}[1]{$(\ref{#1})$}

\newtheorem{theo}{Theorem}

\newtheorem{prop}{Proposition}
\newtheorem{lemm}{Lemma}

\def\N{\mathbb{N}}          

\def\0{{\bf 0}}
\newcommand{\nx}{\eta}

\def\Z{\mathbb{Z}}
\def\R{\mathbb{R}}

\newcommand{\RR}{{\cal B}}
\newcommand{\Rmax}{{R_{\rm max}}}
\newcommand{\Dmax}{{D_{\rm max}}}

\newcommand{\toone}{\stackrel{{L^1}}{\longrightarrow}}

\newcommand{\eqco}{\setcounter{equation}{0}}
\newcommand{\thco}{\setcounter{theo}{0}}
\newcommand{\prco}{\setcounter{prop}{0}}
\newcommand{\laco}{\setcounter{lemm}{0}}
\newcommand{\coco}{\setcounter{coro}{0}}
\newcommand{\cjco}{\setcounter{conj}{0}}

\newcommand{\deco}{\setcounter{defn}{0}}
\newcommand{\allco}{\eqco  \thco \prco \laco \coco \cjco \deco}
\newcommand{\qed}{\rule[-1mm]{3mm}{3mm}}

\newcommand{\Po}{{\cal P}}
\newcommand{\Q}{{\cal Q}}

\newcommand{\Var}{{\rm Var}}
\newcommand{\var}{{\rm Var}}

\newcommand{\diam}{{\rm diam}}

\newcommand{\X}{{\cal X}}
\newcommand{\F}{{\cal F}}

\newcommand{\NN}{{\cal N}}

\newcommand{\eps}{\varepsilon}
\def\bdm{\begin{displaymath}}
\newcommand{\edm}{\end{displaymath}}
\def\benu{\begin{enumerate}}
\def\eenu{\end{enumerate}}
\def\beqn{\begin{equation}}
\def\eeqn{\end{equation}}
\def\be{\begin{equation}}
\def\ee{\end{equation}}
\def\bea{\begin{eqnarray}}
\def\eea{\end{eqnarray}}

\newcommand{\bean}{\begin{eqnarray*}}
\newcommand{\eean}{\end{eqnarray*}}
\newcommand{\bear}{\begin{eqnarray}}
\newcommand{\eear}{\end{eqnarray}}

\newcommand{\tQ}{\tilde{Q}}
\renewcommand{\bar}{\overline}
\renewcommand{\epsilon}{\varepsilon}

\def\hxi{\hat{\xi}}
\def\txi{\tilde{\xi}}
\def\hH{\hat{H}}
\def\hD{\hat{D}}

\def\hnx{\hat{\eta}}
\def\tnx{\tilde{\eta}}
\def\tD{\tilde{D}}
\def\tT{\tilde{T}}

\begin{document}

\title{\bf Growth and roughness of the interface for ballistic deposition}

\author{Mathew D. Penrose  
\\
{\normalsize{\em University of Bath } }
}
\date{August 2006}

\maketitle



\begin{abstract}
In ballistic deposition (BD), $(d+1)$-dimensional
particles fall sequentially at random
towards an  initially flat, large but bounded
 $d$-dimensional surface,
and each particle sticks to the first point of contact.
For both lattice and continuum BD,
 a law of large numbers
in the thermodynamic limit establishes convergence of
the mean height and surface width of the interface to  constants $h(t)$
and $w(t)$, respectively, depending on time $t$. 
We show that $h(t)$ is asymptotically linear in $t$,
while $w(t)$ grows at least logarithmically in $t$ when $d=1$. 
We also give duality results saying that the height above the origin for
deposition  onto an initially flat 
surface is  equidistributed with the maximum height
 for deposition onto  a surface growing
from a single site.

\end{abstract}

\section{Introduction}
\lbl{secintro}
\allco

Scientific interest in
 growth processes associated with the deposition
of particles on surfaces is considerable; see
 Barab\'asi and Stanley \cite{BS}, Cumberland
and  Crawford \cite{CC}, 
 Vicsek \cite{Vik}.
One family of
 deposition models involves
  $(d+1)$-dimensional particles which  rain down sequentially
at random onto a $d$-dimensional substrate (surface); 
when a particle arrives on the existing agglomeration of
deposited particles, it  sticks to the first particle it contacts,
which may result in  lateral growth and `overhangs' 
(if it does not contact any previous particle, it sticks to 
the substrate).
This is known as  {\em ballistic deposition} (BD),
and one reason for studying it is as a more tractable
modification of {\em diffusion limited aggregation}
(see Atar {\em et al.} \cite{AAK}).

The physical sciences literature concerned with ballistic deposition
 in both  lattice and continuum settings  is extensive (see \cite{BS} for
an overview).  As well as numerous simulation studies 
 dating back to Vold \cite{Vold},
this  literature contains analysis by means of scaling theory 
(Family \cite{Fa}, Family and Vicsek \cite{FV},
Kardar {\em et al.} \cite{KPZ}). However, these
arguments have not been made rigorous from the
 mathematician's point of view; see, e.g., page 56 of \cite{BS}.
The rigorous mathematical literature is much less extensive,
but see \cite{AAK,PY2,PY3,Sepp}.
 The present article builds on the approach of Penrose
and Yukich \cite{PY2,PY3}. 

In all versions of BD considered here, the substrate is
the set $\R^d \times \{0\}$, identified with
$\R^d$, or is some sub-region thereof (denoted $Q$).
Thus the substrate  is assumed initially flat.  All particles
are incompressible $(d+1)$-dimensional solids (typically balls or cubes).
Particles arrive sequentially at random positions in
 $\R^d$. When a particle arrives at a position $x \in \R^d$, its centroid
(or some other specified point in the particle identifying its location)
 slides instantaneously down the ray $\{x\} \times [0,\infty)$,
starting from infinity, until the particle hits a position 
adjacent to 
 either the substrate or a previously deposited
particle, at which point its motion stops and it is permanently
fixed. The difference
between lattice and continuum models is that in the lattice
model the  positions at which particles arrive are restricted
to be in the integer lattice $\Z^d$ (embedded in $\R^d$).

One is interested in the height and width (roughness)
of the interface  consisting of `exposed' particles
that are `visible from above'.
 Loosely speaking, the height and width are interpreted
here as the sample mean and standard deviation of the
heights of exposed particles. 
Let $W_{t,n}$ denote the width 
for deposition onto a $d$-dimensional box of volume $n$,
running for ``time'' (i.e., average number of particles
deposited per unit volume of the box) $t$.
Scaling theory \cite{Fa,FV,KPZ} predicts,
 and subsequent experimental and theoretical
studies seem to confirm 
(for an overview, see \cite{BS}, \cite{Vik})
 that there exist a roughness exponent $\alpha$ and
a growth exponent $\beta$, such that
 the surface width
is governed by the
dynamic scaling relation
\begin{equation} \label{S1}
W_{t,n}  \approx n^{\alpha} f(t/n^{\alpha/\beta} )
\end{equation}
where the scaling function $f$ satisfies
$f(x) \propto  x^{\beta},$ 
for $ x \ll 1,$
and $f(x)\approx  C$  for $x \gg 1$.
If the scaling theory is correct, it implies that 
\bean
W_{t,n } \propto
t^\beta
 ~~~~~ t \ll n ;
~~~~~~~~~~~~~
W_{t,n}  
\propto  n^{\alpha},
~~~~~ t \gg n.
\eean
Physicists believe 
(see e.g. p. 56 of \cite{BS}), 
 that  dimension $d = 1$ one has
the `exact' values
$\alpha = 1/2$ and $\beta = 1/3$ (higher dimensional values of
$\alpha$ and $\beta$ are not known). Validating the above
theory  rigorously is a challenge for mathematicians;
in the present work we take some steps in the direction
of rigorously analyzing the regime with $t \ll n$. 

In this paper we consider both lattice and continuum BD models.
 As in \cite{PY3} we first take a thermodynamic limit
(deposition onto an infinite surface). 
In this limit we obtain expressions for the 
limiting height and width, as a function of `time'
(the mean number of particles deposited per unit area).
The next step is to examine the growth of these
functions with time, and in  doing so we go
beyond the continuum analysis in \cite{PY3}.
We show that the limiting height grows asymptotically linearly with time and
that in one dimension, the limiting width grows at least logarithmically. 

The limits are taken in the opposite  order in
\cite{AAK}, while the time-parameter and  the dimensions of the surface
are simultaneously re-scaled in Sepp\"al\"ainen \cite{Sepp}.
Gravner {\em et al.} \cite{GTWjsp,GTW} provide detailed results on a 
 one dimensional discrete-time growth model sharing some
features with BD, which are consistent with the belief that $\beta =1/3$.

\section{Lattice ballistic deposition}
\lbl{seclat}

We first consider a class of lattice ballistic deposition models,
in which all particles are assumed identical.
Let $\0$ denote the origin in $\Z^d$.
Specify a {\em displacement function} $D: \Z^d \mapsto [-\infty, \infty)$
	with the properties that (i) $D(\0) =1 $,
 and (ii) the set  
$\NN := \{ x \in \Z^d: D(x) \neq -\infty\}$ is finite but has at least
two elements (one of which is the origin).
For $x \in \Z^d$ let $\NN_x:= \{x+y: y \in \NN\}$ 
and let $\NN_x^* := \{x-y: y \in \NN\}$.
The set $\NN$ is a `neighbourhood' of the origin
and $\NN_x$ is the corresponding neighbourhood
of $x$. The idea of a displacement function is that
if a particle arrives at $y \in \NN_x$, then it cannot
slide down the ray $y\times [0,\infty)$ below the position
at height $D(y-x)$.

The substrate is represented by a finite subset $Q$ of $\Z^d$ with $|Q|$
elements.  At each site $x \in Q$,
particles arrive at times forming a homogeneous Poisson process of
unit rate, independently of other sites.
We consider two alternative measures of the 
height of the interface at site $x$, the {\em last-arrival} height
$\xi_{t,Q}(x)$   and
the  {\em next-arrival} height 
 $\nx_{t,Q}(x)$. 
The latter is defined in terms of the former
 by 
\bea
\nx_{t,Q}(x) : =
 \max \{ \xi_{t,Q}(y) + D(x-  y): y \in \Z^d \},
\lbl{nxdef}
\eea
where we say  $-\infty + x := -\infty$ for $x \in \R$, so
in fact
$\nx_{t,Q}(x)  = \max \{ \xi_{t,Q}(y) + D(x-  y): y \in \NN_x^*\}$.

The evolution of 
 $\xi_{t,Q}(\cdot)$ proceeds as follows. 
Assume
  $\xi_{0,Q}(z)=0$ for all $z \in \Z^d$, and
as a function of $t$, assume $\xi_{t,Q}(z)$  is right-continuous
and  piecewise constant with jumps only
at the arrival times of the Poisson process of arrivals
at site $z$.  
If  a particle arrives at site $z$ at time $t$,
then the (last-arrival) height
at site $z$ is updated to the next-arrival height immediately before
time $t$, i.e. we set
$$
\xi_{t,Q}(x) = \nx_{t-,Q}(x) := \lim_{s \uparrow t} \eta_{s,Q} (x),
$$
while the last-arrival  heights at other sites remains
unchanged, i.e., $\xi_{t,Q}(y) = \xi_{t-,Q}(y) $ for $y \neq z$. 

Special cases include  the so-called  nearest-neighbour (NN)
and next-nearest neighbour (NNN) models  
\cite{BS}. In the NN model, 
one takes $\NN = \{z \in \Z^d: \|z\|_1\} \leq 1\}$
(i.e., the  origin together with its
 lattice neighbours),
 and the displacement function $D$   
is given by 
 $D(x) =0$ for $x \in \NN \setminus \0$;
this is the version of ballistic deposition considered in
  \cite{Sepp}. In the NNN
model, one takes 
$\NN= \{z \in \Z^d : \|z\|_\infty \leq 1 \}$
(i.e., diagonal neighbours are included)
and takes $D(x)=1$ for all $x \in \NN$; this is the version considered
in \cite{AAK}.

Define the {\em mean height functional} $\bar{\xi}_{t,Q}$
 and {\em width functional}  $W_{t,Q}$
 to be the sample mean and sample
variance, respectively, of the heights at time $t$, i.e., set 
\bea
\bar{\xi}_{t,Q} := |Q|^{-1} \sum_{z \in \Q} \xi_{t,Q}(z); 
~~~~~~
~~~~~~
W_{t,Q}:= 
\sqrt{
|Q|^{-1} \sum_{x \in Q} 
(\xi_{t,Q}(x) - \overline{\xi}_{t,Q})^{2} 
}
.
 \lbl{020723}
\eea
The mean height  $\bar{\xi}_{t,Q}$ is a measure  of the average
amount of empty space under the surface, 
while  the surface width functional 
 $W_{t,Q}$
is a measure of the
roughness of the interface.

We consider a particular limiting regime. First we take
$Q$ to be large with $t$ fixed (deposition for a finite time
onto a very large surface), and then we would like to take
the large-time limit.

The large-$Q$ limit of $\bar{\xi}_{t,Q}$ and $W_{t,Q}$
is best described in terms of  deposition onto  an 
`infinite substrate' represented by the whole of $\Z^d$.
Let $\xi_t(z)$ be the height above site $z$ of this infinite interface
at time $t \geq 0$.
Assume again that $\xi_t(z) =0$ for all $z \in \Z^d$, but
now assume particles arrive as independent
Poisson processes at {\em all} sites in $\Z^d$.
With this as the only difference from the description of
$\xi_{t,Q}(\cdot)$, let the updating rules 
for $\xi_t(\cdot)$ be just the same as before.
Also, define the next-arrival height $\eta_t(x)$ in an
analogous manner to (\ref{nxdef}).
For the infinite interface process
we need to check that no `explosions' occur;
our first result does this and more.
\begin{prop}
 \lbl{th020722a}
  For all $t \in (0, \infty)$, the
values of $\xi_t(\0)$ and $\nx_t(\0)$ are
almost surely finite, and 
for all $k \in \N$, it is the case that
$E[(\xi_t(\0))^k] = O(t^k)$ 
and
$E[(\eta_t(\0))^k] = O(t^k)$ 
as $t \to \infty$.
\end{prop}

We can interpret $\xi_t( \0 ) $ as
 the  height of a `typical' point in the infinite interface.
The next result is
a thermodynamic limit and shows that
 in the large $n$ limit,
 the height
functional and squared width functional
converge
 to the mean and 
variance, respectively, of the `typical' height
$\xi_t(\0)$
(analogously to results in \cite{PY3} for continuum BD).

Let $(Q_n,n \geq 1)$ be a sequence of finite
subsets of $\Z^d$.
Let $\partial Q_n$ denote the set of boundary
sites in $Q_n$, and set
 $\diam(Q_n) := \max_{x,y \in Q_n}|x-y|$.
 Assume that
\bea
\0 \in Q_n 
{\rm ~for~ all~} n;
\lbl{Qn1}
\\
\liminf (Q_n) = \Z^d , ~~i.e.~
\cup_{n=1}^\infty \cap_{m=n}^\infty Q_m = \Z^d; \lbl{Qn2}
\\
 |\partial Q_n|/|Q_n| \to 0
 {\rm~~ as ~~} n \to \infty,
 \lbl{Qn3}
\\
|Q_{n+1}| > |Q_n| 
{\rm ~for~ all~} n; \lbl{Qn4}\\
\diam(Q_n)/|Q_n| \to 0
 {\rm~~ as ~~} n \to \infty.
\lbl{Qn5}
\eea
For example,
$Q_n$ could be a lattice box of side $n$ centred at the origin.

\begin{prop}
\label{th020723}
 For all $t \in (0,\infty)$,
\bea
\bar{\xi}_{t,Q_n}  \toone 
E[\xi_t(\0)] := h(t)
{\rm ~~~~as~~} n \to \infty;
\lbl{020809a}
\eea
\bea 
W^2_{t,Q_n} \toone 
 {\rm Var}[\xi_t(\0) ] := w^2(t)
{\rm ~~~~as~~} n \to \infty.
\lbl{020809b}
\eea
\end{prop}
It is likely that with the Poisson arrivals processes
at all sites in $\Z^d$ all put on the same 
probability space, \eq{020809a} and \eq{020809b}
also hold with almost sure (a.s.) convergence; later we shall give
analogous results for continuum BD with a.s. convergence.

It is possible to give central limit theorems associated with
the above laws of large numbers. For continuum BD, such results
have been given in \cite{PY3}. Similar arguments apply in the
lattice case under consideration here, where one can use the
general lattice central limit theorem
 of Penrose (\cite{P2}, Theorem  3.1).
  
It is of great interest to estimate the limiting  constants
$h(t)$ and $w(t)$ in
Proposition \ref{th020723}, and especially,
to understand the growth of $h(t)$ and $w(t)$
as $t$ becomes large.
Our main results are concerned with this.
 Heuristically, one expects  
$h(t)$ to grow linearly in $t$ since the
expected height should vary directly with the deposition
intensity $t$. The next result demonstrates  this,
and  more.

\begin{theo}
\lbl{th020806}
There is a constant $\rho_1 \in (0,\infty)$ such that
$t^{-1}h(t) \to \rho_1$ as $t \to \infty$.
Moreover, for any $p \in [1,\infty) $ we have the
 $L^p$ convergence 
\bea
t^{-1} \xi_t(\0) 
\to \rho_1   
{\rm ~~~~  as ~} t\to \infty.
\label{020809c}
\eea
\end{theo}

In the  special case of the NN model,  it can be deduced from Theorem 1 of
Sepp\"al\"ainen \cite{Sepp} that \eq{020809c}
holds with almost sure   convergence.
Our approach is somewhat different from that of \cite{Sepp}, 
and is needed for subsequent results.

 As for the width, the scaling theory mentioned in Section \ref{secintro}
predicts that $w^2(t)  = \Theta(t^{2 \beta})$.
Even without scaling theory,
 one expects at least that $w^2(t)= O(t)$ on the basis of 
 simulations of Zabolitzky and Stauffer (\cite{ZS}, p. 1529),
 and also on the following heuristic grounds.
If $\NN =\{\0\}$ then the heights $\xi_{t,n}(x),x \in Q_n$ 
are independent Poisson variables so that $w^2(t) =t$.
If $\NN \neq \{\0\}$ so as to give non-trivial interactions,
these interactions should have a `smoothing' effect
so that $w^2(t)$ should not be any bigger than in
the case $\NN= \{\0\}$.

Rigorous analysis of the large-$t$ behaviour of $w(t)$ appears
to be difficult: the following result makes a start.

\begin{theo}
\label{th020806a}
For any $d$, it is the case that 
\bea
\liminf_{t \to \infty} w^2(t) >0,
\label{fplower0}
\eea
and if $d =1$ and $\NN$ is a lattice interval, then 
\bea
\liminf_{t \to \infty} ( w^2(t)/ \log t) >0.
\label{fplower}
\eea
\end{theo}

The proof of Theorem \ref{th020806a} uses 
a {\em duality} relation between the ballistic 
deposition process and a {\em dual BD} process,
 denoted $(\hxi_t(x), x \in \Z^d)_{t \geq 0}$ with
corresponding next-arrival process $(\hnx_t(x), x \in \Z^d)_{t \geq 0}$,
    defined in an identical manner to the original
BD process
$\xi_t(x)$ and next-arrival process $\nx_t(x)$,
 except that now one uses the
dual displacement function $\hD$ given by $\hD(x) = D(-x)$, 
$x \in \Z^d$,
and takes as initial configuration a single particle
at height 0 at the origin, i.e.,  one takes 
$$
\hxi_0(x) = \left\{ \begin{array}{lr}
0 & \mbox{ if } x =\0 \\ 
 -\infty & \mbox{ otherwise.}
\end{array}
\right.    
$$
We also define a further process $(\tnx_t(x), x \in \Z^d)_{t \geq 0}$,
by
$$
\tnx_t(x) = \hnx_{(t-T)^+}(x)
$$
where
$T$ is exponentially distributed with mean 1, independent of
the process $(\hnx_s(x), x \in \Z^d)_{s \geq 0}$.
In other words, the process $(\tnx_t(\cdot),t \geq 0)$ is
obtained by waiting an exponentially distributed
amount of time before `kicking off' the dual next-arrival
process $\hnx$.  We shall refer to  $\tnx$ as
the {\em delayed dual BD process}.

\begin{theo}
\lbl{th020807}
Let $t>0$.
Then the  distribution of $\nx_t(\0)$ is the same as that
of the maximum depth $\sup_{z \in\Z^d} \hnx_t(z)$  in the dual BD process.
The distribution of $\xi_t(\0)$ is the same as that of 
the maximum depth
$\sup_{z \in \Z^d}
\tnx_t(z)$ in the delayed dual BD process.
\end{theo}
Theorem \ref{th020807} shows that $w^2(t)= \Var( \sup_{z \in \Z^d} \tnx_t(z))$.
The dual BD process is a random interface growing from a single seed; in
this it resembles the classical first passage percolation model.
Our proof of \eq{fplower} uses ideas from the proof by
Pemantle and Peres \cite{PP} 
 of an analogous logarithmic lower bound 
 for  the variance of first passage times.
Further progress in estimating the growth rate of the
variance for first passage times has proved elusive (see e.g. \cite{BKS,Kes});
by analogy, the same could well be true in the case of $w^2(t)$.

\section{Continuum ballistic deposition}
\allco

We  consider a  continuum ballistic deposition model,
defined as follows.  The substrate $\tQ$ is a
Borel-measurable region
of $\R^d$ (for example, a cube of side $n$) and 
$|\tQ|$ denotes the Lebesgue measure of $\tQ$.
Particles are assumed to be $(d+1)$-dimensional Euclidean 
balls of possibly random independent identically distributed  radii
which are uniformly bounded by a finite constant  $\Rmax$
(in fact the results presented here could easily accommodate a
more general class of random shapes, such as convex
shapes of  uniformly bounded diameter). 
Let $F$ denote the common distribution function of the radii of
incoming  particles; assume that $F(0)=0$ and $F(\Rmax)=1$.

Each incoming particle arrives perpendicularly to the  substrate
$\tQ \times \{0\}$ 
and  sticks  to the first previous particle it encounters,
or to the substrate if it does not encounter any previous particle.
In other words, its motion stops when it encounters a previous
particle or the substrate, and 
remains stationary thereafter. 

For simplicity we assume $\tQ$ is given by one of the sets
$\tQ_n (n \geq 1)$  defined by 
$$
\tQ_n := Q_n \oplus [0,1)^d = \{x+y: x
\in Q_n, y \in [0,1)^d \},
$$
with the sequence of finite sets $Q_n \subset \Z^d$ assumed to satisfy the
conditions \eq{Qn1} - \eq{Qn5}. 

Let $\Po $ be a homogeneous Poisson point process
of unit intensity in $\R^d \times [0,\infty)$,
 each point carrying a mark with distribution $F$.
Denote by ${\cal P}_{t,n}$ the restriction of $\Po$ to
$\tQ_n \times [0,t]$, and 
denote by $\Po_t$ the
restriction of $\Po$ to $\R^d \times [0,t]$.

 Represent points in $\Po$ by $(X,T)$, where $X
\in \R^{d}$ denotes the spatial location (center) of the incoming
particle and $T$ its time of arrival. 
Given $n,t$, the BD process driven by $\Po_{t,n}$
is defined as follows. The spatial locations of incoming
particles are given by the spatial locations of the points of $\Po_{t,n}$,
the order in which they arrive is determined by the time-coordinates
of these points (i.e., they arrive in order of increasing time-coordinate),
and their radii are given by the marks the points carry.

\begin{figure}[htbp]
\begin{center}
\includegraphics[width=10cm]{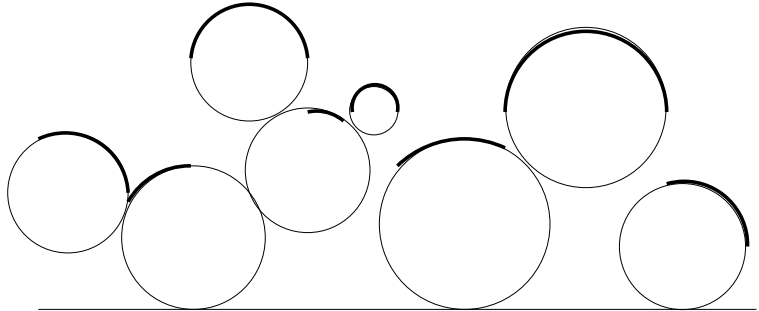}
\end{center}
\caption{The mean value of the function whose graph
is given 
 by the bold arcs is
 $\bar{H}_{t,n}$, and its variance is $W^2_{t,n}$.
The horizontal line represents $\tQ_n \times \{0\}$}
\label{fig0808}
\end{figure}

Let  $A_{t,n}$ denote the agglomeration of particles for
the BD process driven by $\Po_{t,n}$, 
together with the substrate $\tQ_n \times \{0\}$  (a subset of $\R^{d+1}$).
For each  $x \in \tQ_n$ let $H_{t,n}(x)$ denote the height 
of the interface above $x$, i.e.,
let
$$
H_{t,n}(x) := \sup\{h: (x,h) \in A_{t,n}\}.
$$
The bold line in Figure \ref{fig0808}  represents the 
graph of the function $H_{t,n}(\cdot)$.

We define the average height $\bar{H}_{t,n}$  and width $V_{t,n}$
 of the interface at time $t$
as follows.
We set $\bar{H}_{t,n}$
 to be the mean of the function $H_{t,n}(x), x \in Q_n$,
and $V_{t,n}$ to be the root-mean-square deviation of this
function from $\bar{H}_{t,n}$.  That is, we set
\bea
\bar{H}_{t,n} := |Q_n|^{-1} \int_{Q_n} H_{t,n}(x) dx ;  
\lbl{Htndef}
\\
V_{t,n} : = \sqrt{|Q_n|^{-1}  \int_{Q_n} (H_{t,n}(x) - \bar{H}_{t,n})^2 dx}.  
\lbl{Vtn0}
\eea
 These definitions of height and width  are slightly different from those
used in Penrose and Yukich \cite{PY3} but no less natural.

The following result gives meaning to the
height of the interface above 
 an {\em infinite}  substrate.

\begin{prop} \lbl{LLN1} 
 For all $t \geq 0 $, $x \in \R^d$ the limit
\bea
H_t(x) = \lim_{n \to \infty} H_{t,n}(x)
\lbl{Htdef}
\eea
exists almost surely and is almost surely finite. Also,
the distribution of $H_t(x)$ does not depend on $x$.
Moreover, for all $k \in \N$ it is the case that
$E[(H_t(\0))^k] = O(t^k)$ as  $t \to \infty$.
\end{prop}

Thus the infinite-substrate height function is given  by $H_t(x),x \in \R^d$.
The following thermodynamic limits  
are variants of results in \cite{PY3}.
They are continuum analogues to Proposition 
\ref{th020723}.  In particular, they show that the
variance
of the random variable $H_t(\0)$ (height above a typical point of
the infinite substrate)
is the large-$n$ limit  of the width functionals $V^2_{t,n}$ (sample variance
of the heights above a bounded substrate).
\begin{prop}
\lbl{propctstherm}
  For all $t \in (0, \infty)$, we have the 
almost sure convergence
 \bea
 \lim_{n \to \infty} (\bar{H}_{t,n})  = E [ H_t(\0) ] ; 
\lbl{ctslln1}
\\
\lim_{n \to \infty} V^2_{t,n}
=
\var [ H_t(\0) ]. 
\lbl{ctslln2}
\eea
\end{prop}
Now we give a continuum analogue to Theorem \ref{th020806}
\begin{theo}
\label{ctsLLN}
There is a constant $\rho_2 \in (0,\infty)$ such that
$t^{-1}E[H_t(\0)] \to \rho_2$ as $t \to \infty$.
Moreover,  for any $p \geq 1$,
\bea
t^{-1} H_t(\0) 
\to \rho_2  {\rm ~~ in ~~ } L^p{\rm ~~~~  as ~} t\to \infty.
\label{0621b}
\eea
\end{theo}
The next result is the continuum analogue to Theorem.
\ref{th020806a}.

\begin{theo}
\lbl{thctslb}
For any $d$, it is the case that 
\bea
\liminf_{t \to \infty} \var[ H_t(\0) ] > 0 ,
\label{ctslower0}
\eea
and if $d =1$, then 
\bea
\liminf( \Var[H_t(\0)]/ \log t) >0.
\label{ctslower}
\eea
\end{theo}
Also of independent interest is the duality result
given by Proposition \ref{thctsdual} below, which is
the continuum analogue of
Theorem \ref{th020807}.

\section{Proofs for lattice BD}
\allco

For $x \in \Z^d$, let $S_1(x), S_2(x), S_3(x),\ldots$
denote the ordered arrival times of the Poisson process
at site $x$.  The proofs for lattice  BD  
are based on a directed graph representation.
For each $t >0$, define  a directed graph ${\cal G}_t$
and a directed graph $\tilde{\cal G}_t$, 
both with vertex set ${\cal V}_t$ defined by
$$
{\cal V}_t :=
\{(z,S_i(z)): z \in \Z^d, i \in \N, S_i(z) < t\}
 \cup  \{(x,0),(x,t); x \in \Z^d \}.
$$
Informally, each Poisson arrival at site $x \in \Z^d$
before time $t$ is
represented by a point $(x,T) \in \Z^d \times [0,\infty)$.
Two points $(x,T)$ and $(y,U)$ of ${\cal V}_t$
 are joined by a directed
edge in ${\cal G}_t$ from $(x,T)$ to $(y,U)$ if $T<U$ and $y \in \NN_x$. 
They  
 are joined by a directed
edge in $\tilde{\cal G}_t$ from $(x,T)$ to $(y,U)$ if $T<U$ and $y \in \NN^*_x$. 

A {\em path} in ${\cal G}_t$ is a sequence $\pi$ of vertices in ${\cal V}_t$
 denoted $(x_0,T_0),\ldots,(x_n,T_n)$ say,
such that $T_0=0$ and  $T_n = t$ and
for $i=1,2,\ldots, n$ 
there is an edge of ${\cal G}_t$ from
$(x_{i-1},T_{i-1})$ to $(x_i,T_i)$. 
We say the path {\em starts} at $(x_0,0)$ and 
{\em ends} at $(x_n,t)$.
The {\em length} of the path is $n$.
The {\em height} of the path, denoted $h(\pi)$, is
$\sum_{j=1}^n D(x_j-x_{j-1})$.

A path in $\tilde{\cal G}_t$ is defined similarly, except
that now the requirement is that for each $i$
there is an edge of $\tilde{\cal G}_t$ from
$(x_{i-1},T_{i-1})$ to $(x_i,T_i)$, and the 
 (dual) height
 $\hat{h}(\pi)$ of
 a path $\pi = ((x_i,T_i))_{i=0}^n$ in $\tilde{\cal G}_t$ is given by  
$\sum_{j=1}^n \hat{D}(x_j-x_{j-1})$.

The {\em skeleton} of a path $\pi$ in ${\cal G}_t$ is its projection
onto $\Z^d$, i.e., the sequence $(x_0, \ldots, x_n)$. 
Given a sequence $(x_0,\ldots,x_n)\in (\Z^d)^{n+1}$
 with $x_{i}-x_{i-1} \in \NN$ for $1 \leq i \leq n$,
a {\em maximal path} with skeleton $(x_1,\ldots,x_n)$ is
a path $((x_0,T_0),  \ldots, (x_n,T_n))$
in the graph ${\cal G}_t$, 
with $T_0 =0$ and $T_n = t$, and
with the property that for $1 \leq i \leq n$,
 there are no Poisson arrivals at $x_n$ in the time-interval
$(T_{i-1},T_i)$.
\begin{lemm}
\label{le020806}
Let $t >0$ and $z \in \Z^d$. 
Suppose the maximum length of all paths in ${\cal V}_t$ ending at $(z,t)$
is finite. Then $\nx_{t-}(z) := \lim_{s \uparrow t} \nx_s(z)$
is given by
\bea
\nx_{t-}(z) =
\sup \{ h(\pi): \pi ~ \mbox {\rm a path in } {\cal G}_t ~ \mbox{\rm ending at }
 (z,t) \},
\lbl{graph1}
\eea
and $\hnx_{t-}(z) := \lim_{s \uparrow t} \hnx_s(z)$
is given by
\bea
\hnx_{t-}(z) = 
\sup \{ \hat{h}(\pi): \pi ~ \mbox{\rm a path in } \tilde{\cal G}_t ~ 
\mbox{\rm starting at }
(\0,0) ~
\mbox{\rm and ending at }
 (z,t) \},
\lbl{graph2}
\eea
with the convention that the maximum of the empty set is $-\infty$.
\end{lemm}
{\em Proof.}
If $\pi$ is a finite path ending at $(z,t)$,
then we assert  that 
$\eta_{t-}(z) \geq h(\pi)$. Indeed, if $\pi = (x_i,T_i)_{i=0}^n$ 
then by monotonicity of the processes $\xi_t$ and $\nx_t$,
for each $i$ with $1 \leq i <n$ we have 
$$
\xi_{T_i}(x_i) 
\geq \nx_{T_i-}(x_i) \geq \xi_{T_{i-1}}(x_{i-1}) + D(x_i-x_{i-1}) 
$$
so that
\bean
\eta_{t-}(z) \geq \xi_{T_{n-1}}(x_{n-1}) + D( x_n - x_{n-1})
\\
\geq \sum_{i=1}^n D(x_i - x_{i-1}) = h(\pi).
\eean

Conversely, there is at least one path of height  at least $\nx_{t-}(z)$
that ends at $(z,t)$.  This is proved by induction on the maximum
  length of paths ending at $(z,t)$;
it is clearly true when this maximum path length is 1;
  suppose  it is true when the maximum path length is  in the range
$\{1,2,\ldots,  k\}$.  Now suppose 
 the maximum path length  is $k+1$.
Then 
\bean
\nx_{t-}(z) = \max
\{ \xi_{t-}(y) + D(z-y)
: y \in \NN^*_z
 \} \\
\eean
so that for some $ y^* \in \NN_z^*$ (i.e., with $z \in \NN_{y^*}$) we obtain
\bean
\nx_{t-}(z) 
 =
\xi_t(y^*) + D(z-y^*),  
\eean
and if
 $L$ denotes  the last Poisson arrival time before time $t$ at
site $y^*$, we have
\bean
\nx_{t-}(z) 
= \xi_L(y^*) + D(z-y^*) 
\\
= \nx_{L-}(y^*) + D(z-y^*). 
\eean
 Each path in ${\cal G}_L$ ending  at  $(y^*,L)$ 
has length at most $k$, since otherwise there would be a path 
through $(y^*,L)$ to
$(z,t)$ of length greater than $k+1$. Hence by the inductive
hypothesis, there is a path in ${\cal G}_L$
ending at $(y^*, L)$, of height at least
$\nx_{L-}(y^*)$, and hence by appending $(z,t)$ to this sequence
one obtains a path in ${\cal G}_t$, ending at
 $(z,t)$, of height at least $\nx_{L-}(y^*) + D(z-y^*)$.
This completes the induction and hence the proof \eq{graph1}.
The proof of \eq{graph2} is similar.
$\qed$ \\  

\noindent
{\em Proof of Proposition \ref{th020722a}.}
Since $\xi_t(\0) \leq \eta_{t}(\0)$, it suffices 
to prove that $\nx_t(\0)$ is a.s. finite and
$E[\nx_t(\0)^k]= O(t^k)$ as $t \to \infty$.
Note that $\nx_t (\0) = \nx_{t-}(\0)$ a.s. 
Also, let $\Dmax : = \max_{z \in \Z^d} D(z)$ denote
the maximum value of the displacement function.
By our assumptions on this function, we have
$1 \leq \Dmax < \infty$.

By Lemma \ref{le020806},
if $\nx_{t-}(\0) \geq m$
then there must be a path in ${\cal G}_t$
 of height at least $m$, and hence of length
at least $m/\Dmax$, that ends at $(\0,t)$. 
Hence, if $\nx_{t-}(\0) \geq m$ then there is a path 
in ${\cal G}_t$ of length at least $m/\Dmax$ that
is maximal for its skeleton and ends at $(\0,t)$. 

For any given sequence $(x_0,x_1,\ldots,x_n)$ with
$x_{i} -x_{i-1} \in \NN$ for $i=1, \ldots, n$,
the probability that there exists a maximal path in
${\cal G}_t$ with skeleton 
$(x_0,x_1, \ldots, x_n)$
and with all arrival times less
than $t$, equals than the probability that the sum
of $n-1$ independent $\exp(1)$ variables is less than $t$,
which is the same as the probability that ${\rm Po}(t) \geq n-1$,
where ${\rm Po}(t)$ denotes a Poisson variable with mean $t$.
Hence, for any $y \geq t {\rm e}^2+1$,
by e.g. Lemma 1.2 of \cite{Pen},
\bea
P[\nx_{t}(\0) \geq y \Dmax] 
=
P[\nx_{t-}(\0) \geq y \Dmax] 
\leq |\NN|^y P[ {\rm Po}(t) \geq y-1]
\nonumber \\
\leq |\NN|^y \exp( - ((y-1)/2) \log ((y-1)/t)) 
\lbl{060703}
\eea
which tends to zero as $y \to \infty$.
Hence $\nx_t(\0)$ is almost surely finite.
Also, for
$k \in \N$,
\bean
E\left[ \left(\frac{ \nx_t(\0)}{\Dmax} \right)^k \right] 
= \int_0^\infty P \left[ \frac{\eta_t(\0)}{\Dmax}
\geq w^{1/k} \right] dw.  
\eean
Set $c := (2 |\NN|)^2$, split the region of integration into
 $w \leq (ct \leq 1)^k$
and
$w \geq (ct +1)^k$, and change variables to $y=w^{1/k}$
 in the second integral to obtain the estimate
\bean
E\left[ \left(\frac{ \nx_t(\0)}{\Dmax} \right)^k \right] 
\leq (ct+1)^{k} +
 \int_{ct+1}^\infty P \left[  \frac{\eta_t(\0)}{\Dmax} \geq y
\right]
ky^{k-1} dy
\\
\leq (ct+1)^{k} +
k {\rm e}^{1/2}
 \int_{ct+1}^\infty 
 |\NN|^y
 \exp(- (y/2) \log (c)) 
y^{k-1}
dy
\\
\leq (ct+1)^k 
+ k {\rm e}^{1/2}
  \int_{0}^\infty y^{k-1}  \exp(- y \log 2) dy
\eean
and the last integral is finite, so that
 $t^{-k} E[\nx_t(\0)^k]$ is bounded.
$\qed$ \\

\noindent
{\em Proof of Proposition \ref{th020723}.}
The idea is to apply Theorem 3.1 of Penrose \cite{P2}. 
Let $\RR$ denote the collection of all non-empty finite subsets of
$\Z^d$.
Let $X=(X_x, x \in \Z^d)$ be a family of independent homogeneous
Poisson processes of unit intensity. For
$Q \in \RR$, assume  the evolution of $\xi_{t,Q}$ is governed
by the Poisson processes $(X_x,x \in Q)$. Then $(\xi_{t,Q}(x), 
Q \in \RR,  x \in Q)$
is stationary $\RR$-indexed summand in the sense of
section 3 of \cite{P2}.
Also, the proof of Proposition 2.1 shows that 
\bea
\sup \{
 E[ \xi_{t,Q}(x)^4 ]: 
Q \in \RR, x \in Q
\}
 < \infty.
\lbl{mom4}
\eea
For any $\RR$-valued sequence $(B_n, n \geq 1)$
 with $\liminf(B_n) =\Z^d$, we have
as $n \to \infty$ that
  $\xi_{t,B_n}(\0) \to \xi_t(\0)$ almost surely 
(see the proof of Proposition 2.1, or Lemma 5.1 of \cite{P2}),
and hence in  $L^2$  by \eq{mom4}.
Hence, for any such sequence $(B_n, n \geq 1)$
we have the $L^1$ convergence
\bea
\xi_{t,B_n}(\0) \to \xi_t(\0); \lbl{pencon1}
\\
\xi_{t,B_n}(\0)^2 \to \xi_t(\0)^2.
\lbl{pencon2}
\eea
By  \eq{pencon1},
 the first part of Theorem 3.1 of  \cite{P2}  
is applicable to the stationary $\RR$-indexed summand  
$(\xi_{t,Q}(x),Q \in \RR,x \in Q)$, 
so that
eqn (3.3) of  \cite{P2} gives us  
\eq{020809a}.

Moreover, since \eq{pencon1} holds with $L^2$ convergence as well,
by changing $L^1$ estimates to $L^2$ estimates throughout
the proof of eqn (3.3) of \cite{P2}, we may deduce
that result in this case with $L^2$ convergence, i.e. \eq{020809a} holds
with $L^2$ convergence as well (the proof in \cite{P2} in its turn 
uses a multiparameter $L^1$ ergodic theorem quoted from \cite{P0},
but this is also easy to extend to $L^2$ convergence 
in the present setting, using \eq{mom4}.)

To prove \eq{020809b}, first  expand the sum of squares in
\eq{020723} to obtain
\bea
W_{t,Q_n}^2 = |Q_n|^{-1} 
\left(  \sum_{x \in Q_n} \xi_{t,Q_n}(x)^2 \right) - \bar{\xi}_{t,Q_n}^2.  
\lbl{varexpand}
\eea
By \eq{pencon2}, the first part of Theorem 3.1 of  \cite{P2} is
 applicable 
to the stationary $\RR$-indexed summand
$(\xi_{t,Q}(x)^2,Q \in \RR,x \in Q)$, 
and  this shows that the first term in the right hand side of 
\eq{varexpand} converges in $L^1$ to $E[ \xi_t(\0)^2]$.

Since \eq{020809a} holds with $L^2$ convergence, it follows
that the second term in the right hand side of \eq{varexpand}
converges in $L^1$ to
 $(E[ \xi_t(\0)])^2$.
Combining these limiting results in \eq{varexpand}, we obtain
\eq{020809b}.
$\qed $\\

{\em Proof of Theorem \ref{th020807}.}
The idea of the proof  is a form of {\em time-reversal}
of the graphical representation.
Let $\psi_t: \Z^d \times [0,t] \to \Z^d \times[0,t]$ be defined by
$$
\psi_t((z,s)) = (z,t-s).
$$
Let $\hat{\cal V}_t: = \psi_t(  {\cal V}_t)$,
and observe that the  point set $\hat{\cal V}_t$ has the same distribution 
as ${\cal V}_t$ by the properties of Poisson point processes.

Let $\hat{\cal G}_t$ be defined in the same manner as $\tilde{\cal G}_t$
but on the vertex set $\hat{\cal V}_t$ instead of ${\cal V}_t$.
Then each path in ${\cal G}_t$ ending at $(0,t)$
corresponds to a path with the same  height  starting at $(0,0)$  
in $\hat{\cal G}_t$; the correspondence is obtained by reversing the
sequence of vertices  and then applying the mapping $\psi_t$ to each
vertex in the sequence.

Note that $\nx_t(\0)= \nx_{t-}(\0)$ with probability 1.
By \eq{graph1},
$\nx_t(\0)$ is the greatest height of
all paths in ${\cal G}_t$ which end at $(\0,t)$. 
By the correspondence described above, this is precisely the same as
the greatest height of all paths  
in the graph $\hat{\cal G}_t$ which start at $(\0,0)$. 

Since the point processes $\hat{\cal V}_t$ ${\cal V}_t$ have the same
 distribution, it follows that $\nx_t(\0)$ has the same distribution as
 the greatest height of all paths in the  graph $\tilde{\cal G}_t$
 which start at $(\0,0)$. Hence by \eq{graph2},
 $\nx_t(\0)$ has the same distribution as
$\sup_{z \in \Z^d} \hnx_{t-}(z)$ and hence
the same distribution as
$\sup_{z \in \Z^d} \hnx_{t}(z)$.

Let $L$ be the last arrival time at $\0$ before time $t$ (or $L=0$ 
if there are no arrivals at $\0$ before time $t$).
Then $t-L$ has 
the distribution of $\min(T,t)$ where  $T$ is
exponential with mean 1, and hence $L$ has the
distribution of $(t-T)^+$.  Given $L $ with $L>0$, 
the distribution of $\nx_{L-}(\0)$ is the same as that of 
$\sup_{z \in \Z^d} \hnx_{L}(z)$, by the same argument as above
and the Markov property of the time-reversed 
Poisson process.
Whenever $ L >0$  
we have
 $\xi_t(\0)= \nx_{L-}(\0)$, and if $L=0$ then $\xi_t(\0) =0 $ almost surely. 
Combining these observations shows that $\xi_t(\0)$ has
the same distribution as $\sup_{z \in \Z^d} \hnx_{(t-T)^+}(z)$,
as asserted.
$\qed$ \\

From now onwards, we shall assume the delayed dual BD (next-arrival)
 process $\hnx_t(z)$ 
and also the associated last-arrival process $\hxi_t(z)$
is defined in terms of the Poisson arrival times 
$(S_i(x), i \geq 1, x \in \Z^d)$, as follows. For $t < S_1(\0)$
 we put $\hnx_t(z) =0$ and $\hxi_t(z) =-\infty$ for all $z \in \Z^d$.
Then we put 
$$
\hxi_{S_1}(z) :=    
 \left\{ \begin{array}{lr}
0 & \mbox{ if } z =\0 \\ 
 -\infty & \mbox{ otherwise.}
\end{array}
\right.    
$$
and define $\hnx_{S_1(\0)}$ in terms of $\hxi_{S_1(\0)}$ in the usual manner
as given at \eq{nxdef}.
Then we allow the
evolution of $(\hxi_t,\hnx_t)_{t \geq S_1(\0)}$ 
to follow the usual rules driven by the Poisson arrivals
$\{S_i(z): S_i (z) > S_1(\0)$.
Then $\hnx_t(x),x \geq 0$ constructed in this manner
 clearly follows the desired evolution prescribed
in Section \ref{seclat}, with $S_1(\0)$ used as the
initial `kicking off time'. 

For $t \geq 0$, $u \geq 0$, let $D_t$ be the depth 
(i.e., the maximum next-arrival height) of the 
delayed dual BD process
 $\tnx_t$ at time $t$, and let $T(u)$ be first passage time 
to depth $u$ for  the delayed dual BD process, i.e.
  let
\bean
D_t : = \sup_{z \in \Z^d} \{ \tnx_t(z) \}; ~~~
T(u) : = \inf \{t: D_t \geq u\}.
\eean

\begin{lemm}
\lbl{lemmaxlim}
There is a constant $\rho \in (0,\infty)$ such that
\bea
\lim_{u \to \infty}
\frac{T(u)}{u} = \rho, ~~~a.s.
\lbl{060620}
\eea
and 
\bea
\lim_{t \to \infty}
\frac{D_t}{t} = \rho^{-1}, ~~~a.s.
\lbl{060815}
\eea
\end{lemm}
{\em Proof.}
First we verify that $E[T(1)^2]$ is finite. Note that 
$\tnx_t(\0) \geq N_t-1$, where here $N_t $
denotes the number of arrivals at $\0$ up to time $1$.
Hence, 
\bean
P[T(1) > t] \leq P[ \tnx_t(\0) < 1 ] \leq P[N_t < 2] = P[ X+Y >t]
\eean
where here $X$ and $Y$ denote independent exponential random variables
with unit mean, representing the first two inter-arrival times at $\0$.
Hence,
$$
E[ T(1)^2] \leq E[ (X+Y)^2] < \infty.
$$

Next, we  assert that  
 $T(u)$ is 
distributionally subconvolutive, i.e. for $u, v \geq 0$ we have
\bea
F_{T(u+v)} \geq F_{T(u)} * F_{T(v)}.
\lbl{subconv}
\eea
To see this, let $X(u) \in \Z^d$ be chosen 
(in an arbitrary way if there is more than one choice)
so that $\tnx_{T(u)}(X(u)) \geq u$ (by definition such $X(u)$ exists).
Let $T(u) + T^*$ be the time of next Poisson arrival after
$T(u)$ at site $X(u)$, and
 let $(\hxi^*_s,\hnx^*_s)_{s \geq 0}$ be a version of the BD process
with displacement function $\hat{D}$, 
with  initial profile
$$
\hxi^*_0(x)= 
 \left\{ \begin{array}{lr}
0 & \mbox{ if } x =X(u) \\ 
 -\infty & \mbox{ otherwise}
\end{array}
\right.    
$$
and driven by Poisson arrival times $\{S_i^*(x)\}$ given for each $ x \in \Z^d$
 by 
$$
\{S_i^*(x) \}=  \{S_i(x) - T(u)-T^*: S_i(x) > T(u)+T^*\},
$$
where $\{S_i(x)\}$ are the arrival times driving the
original dual BD process $(\hxi_t,\hnx_t)$.  

Let $T^{**}(v)$ be the first time the process 
$(\hnx_s^*)_{s \geq 0}$ achieves
a depth of at least  $v$, i.e. 
$$
T^{**}(v) = \inf\{ t \geq 0: \sup_{x \in \Z^d} 
(\hnx^*_t(x)) \geq v \}.
$$
Then $T^* + T^{**}(v)$ has the same distribution as $T(v)$,
and is independent of $T(u)$. Also, since $\txi_{T(u) +T^*}(X(u)) \geq u$,
 the depth at time
 $T(u) + T^* + T^{**}(v)$ is at least $u+v$, i.e.
 $$
\sup_{x \in \Z^d} (\tnx_{T(u)+ T^* + T^{**}(v)}(x) )
\geq u+v,
$$
 so that $T(u) + T^* + T^{**}(v) \geq T(u+v)$.
Combining these facts gives us (\ref{subconv}).
Since the variables $(T(u), u \geq 0)$ 
 are also
 monotonically increasing in $u$, we can apply 
 the Kesten-Hammersley theorem
(\cite{SW}, page 20) to
obtain  (\ref{060620}), with $0 \leq \rho < \infty$.
Also,
\bea
\frac{D_t}{T(D_t)} \geq
 \frac{D_t}{t} \geq \frac{D_t}{D_t+1} \times \frac{D_t+1}{T(D_t +1)} 
\lbl{060815a}
\eea
and since $D_t \to \infty$ almost surely this
with \eq{060620}
yields \eq{060815} provided $\rho >0$.

If $\rho =0$ then by \eq{060620} and  the second inequality of \eq{060815a}
we would have $D_t/t \to \infty$ almost surely, so that
$E[D_t/t] \to \infty$, and so by Theorem \ref{th020807},
$E[\xi_t(\0)/t ] \to \infty$.  
This contradicts Proposition \ref{th020722a}, and hence
 $\rho >0$ as asserted.
$\qed$ \\

{\em Proof of Theorem \ref{th020806}.}
Set $\rho_1 := \rho^{-1}$, with $\rho$ as given in Lemma \ref{lemmaxlim}. 
By Theorem \ref{th020807}, 
$D_t$ has the same distribution as $\nx_t(\0)$,
 and so by  Proposition \ref{th020722a}, 
for any $p \geq 1$ the $p$th moment of
$(D_t/t) $ is  bounded uniformly in $t$.
By \eq{060815}, $(D_t/t)$ converges almost surely to $\rho_1$,
and by the moment bound 
 the almost sure convergence 
extends to convergence in  $p$th moment for any $p \geq 1$.
Since $D_t$ has the same distribution as $\nx_t(\0)$,
this  convergence in $p$th moment also holds for $\nx_t(\0)$. $\qed$ \\

{\em Proof of Theorem \ref{th020806a}.}
To prove (\ref{fplower0}) consider the
event that (i) there are no arrivals in $\NN \setminus \{\0\}$
between times $t-1$ and $t$, and (ii) there is at least one arrival
at $\0$ in the time-interval $(t-1,t]$.
This event has the same non-zero probability for all $t \geq 1$.
 Conditioned on this event, and on everything 
before time $t-1$, the conditional variance of $\xi_t(\0)$
is the variance of the number of Poisson arrivals at $\0$
in the time-interval $(t-1,t]$, i.e. the variance of
a Poisson variable with unit mean  conditioned to take a value
of at least 1. This variance is a strictly positive constant and
(\ref{fplower0}) follows.

Now suppose $d=1$ and $\NN$ is a lattice interval.
To prove (\ref{fplower}),  we consider the delayed dual BD process
$(\txi_t,\tnx_t)_{t \geq 0}$.  In this process, we denote
by {\em accepted arrival}  an arrival 
 time $T = S_i(z)$ such that $\txi_T(z) > -\infty$.
Enumerate
the set of all accepted arrival times (for all sites)
in increasing order
as $\tau_1,\tau_2,\ldots$.
Let $N_t$ be the number of accepted arrivals up to time $t$, i.e.
$$
N_t := \sup\{n: \tau_n \leq t\}.
$$

Let $I_t$ be the  size of the interface at time $t$,
i.e., the number of sites  $z \in \Z$ with $\tnx_t(z) >-\infty$.
For $n \geq 1$, let $Y_n:= I_{\tau_n}$ be the size of the interface 
 after $n$ accepted arrivals, and set $Y_0:=1$.
Since $(N_t, t \geq 0)$ is a Poisson counting process
with its `clock' running at speed $I_t$,
we have that
\bea
\lim_{t \to \infty}
\frac{N_t}{\int_0^t I_u du} =  1, ~~~ {\rm a.s.}
\label{060620a}
\eea
Next, we assert that there is a constant $\gamma>0$ such that
\bea
 \lim_{t \to \infty}
(I_t/t) =  \gamma, ~~~ 
 {\rm a.s.}
\label{060620d}
\eea
To see this, recall  that we are assuming here  that $\NN$ is a
lattice interval including $\0$ and at least one other element.
It is not hard to see that $I_t$ must also be a lattice interval,
and that both the right and the left endpoint of $I_t$
follow renewal reward processes, where in both cases the
 the inter-arrival times of the underlying renewal process
are  exponentially distributed and
the rewards are uniformly distributed over a lattice interval.
The assertion \eq{060620d} follows by
the Strong
Law of Large Numbers for a renewal reward process.

By (\ref{060620a}) and \eq{060620d} we obtain
\bea
\lim_{t \to \infty} \frac{N_t}{t^2} = \frac{\gamma}{2}, ~~ {\rm a.s.}
\label{060620b}
\eea
Let $M(u)$ be the number of accepted arrivals up to time
$T(u)$.
By (\ref{060620b}) and
 (\ref{060620}),
 as $u \to \infty$, it is the case
 with probability 1 that
\bea
M(u) = N_{T(u)} \sim  (\gamma/2) T(u)^2  \sim (\gamma \rho^2/2) u^2.
\lbl{060620e}
\eea

Let ${\cal F}$ be the $\sigma$-algebra
generated by  the locations in $\Z$ of the sequence of accepted arrivals.
Conditional on ${\cal F}$, the distribution of $T(u)$
is that of the sum of $M(u)$ independent exponentials with the
$j$th exponential having mean $Y_{j-1}^{-1}$. 
Hence,
\bea
\sigma_u^2 := \Var[T(u)|{\cal F}]
 = \sum_{j=1}^{M(u)} Y_{j-1}^{-2}.
\lbl{060620c}
\eea
By definition, $N_{\tau_j}=j$, and hence by (\ref{060620b}), 
$j/\tau_j^2 \to \gamma/2$ so that $\tau_j \sim (2 j/\gamma)^{1/2}$ 
as $j \to \infty$, almost surely.  Since $Y_j = I_{\tau_j}$,
by (\ref{060620d})
we obtain as $j \to \infty$ that
with probability 1,
\bea
Y_j \sim  \gamma \tau_j \sim (2 \gamma j)^{1/2}, 
\label{060621a}
\eea
so that by \eq{060620c} and \eq{060620e},
as $u \to \infty$ we have that
\bea
\sigma_u^2 \sim
 (2 \gamma)^{-1} \log M(u) 
\sim \gamma^{-1} \log u, ~~~~~a.s.
\lbl{sigsim}
\eea

The Berry-Esseen theorem, e.g. as given in theorem 5.4 
of Petrov \cite{Petr} or  as quoted in Chen and Shao \cite{CS},
says that there is a constant $C$ such that
 if $X_1,\ldots X_k$ are independent random 
variables with mean zero and finite third moments and
 $W := \sum_{i=1}^k X_i$ has variance 1, then 
$$
\sup_{x \in \R} \{ | P[W \leq x] -\Phi(x)| \} \leq C \sum_{i=1}^k
E[|X_i|^3], 
$$
where $\Phi$ is the standard normal distribution function.

Let $\tau_0 :=0 $, and for
 $i \geq 1$
let $e_i := \tau_i - \tau_{i-1}$.
As mentioned earlier, conditional on ${\cal F}$ the 
$e_i$ are independent exponentials with $E[e_i|\F] = Y_{i-1}^{-1}$.
Define 
\bea
\theta_u :=  \sum_{i=1}^{M(u)} E[|e_i - E[e_i|{\cal F}]|^3|{\cal F}]
\nonumber \\
=  (12{\rm e}^{-1} -2) \sum_{i=1}^{M(u)} Y_{i-1}^{-3},
\lbl{twelve}
\eea
since
if
 $X$ is exponential with mean $a$ then $E[|X-a|^3] = a^3 (12{\rm e}^{-1} -2)$.

Set $\mu_u : = E[ T(u)|{\cal F}]$.  
By the Berry-Esseen Theorem,
\bea
\sup_{x \in \R} \left|  P \left[ \frac{T(u) - \mu_u}{\sigma_u} \leq x  | {\cal F} \right] - \Phi(x) \right| \leq
 \frac{ C \theta_u}{\sigma_u^{3}}.
\lbl{BE}
\eea
By (\ref{060621a}),
the sum $\sum_{i=1}^\infty Y_{i-1}^{-3}$ converges
almost surely. Hence by  \eq{sigsim} and \eq{twelve},
 we can find $u_0$ such that for $u \geq u_0$ we
have $P[A_u] < 0.01$ where 
$$
A_u := 
 \left\{ 
 \frac{ C \theta_u}{\sigma_u^{3}}
> 
0.01
\right\} \cup \{ \sigma_u <    \sqrt{(\log u)/(2 \gamma)} \}
$$
so that $A_u \in {\cal F}$.
Then for any $y \in \R$,  using \eq{BE} we may deduce that
\bean
P[T(u) \leq y  + 0.2 (\gamma^{-1} \log u)^{1/2} |A_u^{\rm c} ] 
 - P[ T(u) \leq y  |A_u^{\rm c}]
\\
=
P\left[\frac{ T(u) -\mu_u }{\sigma_u}   \leq 
\frac{y - \mu_u}{\sigma_u}
+ \frac{0.2  (\gamma^{-1} \log u)^{1/2}}{ \sigma_u} 
 |A_u^{\rm c} \right] - P
\left[ \frac{ T(u) -\mu_u}{\sigma_u} \leq \frac{y-\mu_u}{\sigma_u} 
|A_u^{\rm c} \right] 
\\
\leq P \left[ 
\frac{ T(u) -\mu_u }{\sigma_u}   \leq 
\frac{y - \mu_u}{\sigma_u}
+ (0.2)\sqrt{2}
|A_u^{\rm c}
\right]
-
 P \left[ 
\frac{ T(u) -\mu_u }{\sigma_u}   \leq 
\frac{y - \mu_u}{\sigma_u}
|A_u^{\rm c}
\right]
\\
\leq \sup_{x \in \R^d} \{ \Phi(x + (0.2)\sqrt{2}) - \Phi (x)\} + 0.02
\leq 0.22,
\eean
and so for  $u \geq u_0$ we have
\bea
P[T(u) \leq y + 0.2 (\gamma^{-1} \log u)^{1/2}  ] 
 - P[ T(u) \leq y ] 
\leq
 0.22 + P[A_u] < 1/4.
\lbl{0624}
\eea

 For $t >0$, let $\nu(t) $ be the median
of the distribution of $D_t$. Then 
\bea
P[ T(\nu(t) ) > t] = P[D_t < \nu(t)] \leq 1/2
\label{020810g}
\eea
so that $P[T(\nu(t)) \leq t ] \geq 1/2$, and hence 
by (\ref{0624}),
\bean
P[T(\nu(t)) \leq  t- 0.2 (\gamma^{-1} \log \nu(t))^{1/2} ] \geq 1/4.
\eean
By \eq{060815}, $(D_t/t) \to \rho^{-1}$ in probability, so that
$(\nu(t)/t) \to \rho^{-1}$ as $t \to \infty$. Hence, 
\bea
\log \nu(t) \sim \log t 
~~ \mbox{\rm as }~
 t \to \infty.
\lbl{060815b}
\eea

If $T(\nu(t)) \leq t - 0.2 (\gamma^{-1} \log \nu(t))^{1/2}$,
then there is a site $z^* \in \Z$ with 
$$
\tnx_{t-0.2(\gamma^{-1} \log \nu(t))^{1/2}}(z^*)
\geq \nu(t).
$$
 If also at least $0.1 (\gamma^{-1} \log t)^{1/2}+1$
  particles arrive at site $z^*$ 
between times $t- 0.2 (\gamma^{-1} \log \nu(t))^{1/2}$ and $t$, then
$D_t \geq \nu(t) + 0.1 (\gamma^{-1} \log t)^{1/2}$.
By \eq{060815b},
the conditional probability of the second of these events, given
 the first,  tends to 1
as $t \to \infty$, 
so that for large $t$,
\bea
P[D_t \geq  \nu(t) + 0.1 (\gamma^{-1} \log t)^{1/2} | 
> 1/8.
\lbl{0624a}
\eea
Moreover, by definition $P[D_t \leq \nu(t)] \geq 1/2$. Combining this
wtih \eq{0624a} yields (\ref{fplower}).  $\qed$ \\

\section{Proofs for continuum BD}
\allco
In this section, for $x \in \R^d$ and $r >0$ we
write $B_r(x)$ for the closed Euclidean ball of radius $r$
centred at $x$.

We introduce a {\em dual continuum BD} process, denoted $\hH_t(x)$,
a    defined in an identical manner to the original
continuum BD process except that now the initial profile
$\hH_0(x), x \in \R^d$ is given by
$$
\hH_0(x) = \left\{ \begin{array}{lr}
0 & \mbox{ if } x =\0 \\ 
 -\infty & \mbox{ otherwise.}
\end{array}
\right.    
$$
In other words, the interface grows from an initial  seed
consisting of a single point
at $(\0,0)$.
Incoming  particles which miss the agglomeration 
have no effect.

We adapt the graphical representation to the continuum.
 Given a locally finite  marked point set $\X \subset \R^d \times [0,\infty)$,
we say a sequence of points $(X_i,T_i)_{i=1}^k$ in $\X$
forms a {\em path in $\X$} if $T_1 < T_2 < \cdots <T_k$, and
 for $1 \leq i <k$, the $d$-dimensional ball centred at $X_i$
 (of radius given by the mark of $(X_i,T_i)$) overlaps 
 the ball centred at $X_{i+1}$
 (of radius given by the mark of $(X_{i+1},T_{i+1})$).

Given such a path,
we refer to the sequence $(X_1,\ldots,X_k)$ as
the {\em skeleton} of the path.
Given also $x \in \R^d$, we say the path {\em starts  near} $x$ if $x$ lies in
the $d$-dimensional ball centred
at $X_1$ (of radius given by the corresponding mark)
and we say the path {\em ends  near} $x$ if $x$ lies in
the $d$-dimensional ball centred
at $X_n$.
 We say the {\em height} at $x$ of the path
is the height of the interface above  $x$  
for the BD agglomeration onto an initially flat surface
 determined by the finite
sequence of incoming $(d+1)$-dimensional balls
 centred at $X_1,\ldots, X_k$ (with radii
given by the corresponding marks).

If the path starts near $\0$, we say the {\em dual height}
of the path at $x$ is the height of the interface above  $x$  
for the BD agglomeration determined by the finite
sequence of incoming $(d+1)$-dimensional balls
 centred at $X_1,\ldots, X_k$ (with radii
given by the corresponding marks), starting initially 
 the initial profile being
given by the function $\Hat{H}_0$ (i.e. by a single point at $(\0,0)$
rather than an initially flat surface).

\begin{lemm}
\lbl{lempath}
Let $t \geq 0$, $x \in \R^d$.
Then  with probability 1:

 (i) There are a.s. only finitely many paths
in $\Po_{t} $
which end near $x$; 

(ii) the value of 
$H_t(x)$ is the maximum height at $x$ of
all such paths (or is zero if no such path exists);

  (iii) $H_{t,n}(x)$ is
 the maximum height at $x$ of
all paths in $\Po_t$ which end near $x$
  for which the skeleton is contained in $\tilde{Q}_n$,
(or is zero if no such path exists), and

(iv) $\hat{H}_t(x)$ is the maximum 
{\em dual} height at $x$ of all
 paths which start near $\0$ and end near $x$ 
(or is $-\infty$ if no such path exists).
\end{lemm}
{\em Proof.}
Part (i) follows from 
 e.g. corollary 3.1 of \cite{P0}.

 In the continuum BD process,
 there exists a sequence of particles in the agglomeration,
each particle touching the next one in the sequence,
 leading from the substrate to a  particle arriving at
 time at most $t$ and with $(x, H_t(x))$ on
its surface in the agglomeration;
considering only the particles in 
this sequence, we have a path which ends near $x$  and has
 height $H_t(x)$ at $x$.
 Hence the height $H_t(x)$ is at most the maximum height
at $x$ of all such paths.

 On the other hand, inserting extra
points into a given $(d+1)$-dimensional marked
point process 
 cannot decrease the height over $x$ of the interface 
of the corresponding  BD agglomeration, 
 so for any path in $\Po_t$ ending near $x$,
$H_t(x)$ is at least the height at $x$ of
the path.
This demonstrates part (ii).

The argument for part (iii) is the same as for part (ii),
except that now one ignores particles with spatial location
outside $\tilde{Q}_n$.

The argument for part (iv) is the same as for part (ii).
$\qed$ \\

{\em Proof of Proposition \ref{LLN1}.}
Let $x \in \R^d$.
By \eq{Qn2} and part (i) of Lemma \ref{lempath}, there
exists an almost surely finite random $n_0$ such that
for $n \geq n_0$, every  path in $\Po_t$ which ends near $x$
has its skeleton contained in $\tilde{Q}_n$. Then by parts (ii) 
and (ii) of Lemma \ref{lempath}, we have $H_t(x) = H_{t,n}(x)
$ for $n \geq n_0$. This demonstrates the
 first part of Proposition \ref{LLN1}.

Since $\Po_t$ is distributionally invariant under
spatial translations (i.e., translations of $\R^{d} \times \R$ 
leaving the time-coordinate
unchanged), the distribution of $H_t(x):= \lim_{n \to \infty} H_{t,n}(x)$
does not depend on $x$.
This demonstrates the second part of Proposition \ref{LLN1}.

We prove the last part only in the case where $\Rmax \leq 1/2$;
the more general case can then be deduced by some simple scaling
arguments which we omit.

We couple our continuum BD process to a certain NNN lattice BD model,
as defined in Section \ref{seclat}.
Partition $\R^d$ into half-open unit cubes, and
 for $x \in \R^d$ let $Q(x)$
be the the cube in the partition that contains $x$.
Let $(\xi_t, \nx_t)_{t \geq 0}$
  be the coupled NNN lattice BD model, in which the arrival times
at site  $z \in \Z^d$ are given by the time-coordinates of the points
of $\Po$ in $Q(z) \times (0,\infty)$. 
By the assumption that $2\Rmax \leq 1$, it is not
hard to see that $H_t(x) \leq \eta_t(z(x))$. 
We can then use Proposition \ref{th020722a} to
deduce that  $E[(H_t(\0))^k] = O(t^k)$. 
$\qed$ \\

{\em Proof of Proposition \ref{propctstherm}.}
By Lemma \ref{lempath} and distributional invariance of $\Po_t$
under spatial translations, there is a function $q(t), t \geq 0$
such that $q(t) \to 0$ as $t \to \infty$ and
\bea
P[H_t(x) \neq H_{t,n}(x) ] \leq q(  {\rm dist}(x, \partial Q_n)). 
\lbl{expest}
\eea
Also, by the argument in the proof of the last part of Proposition
\ref{LLN1}, for any $k \geq 1$ we have 
\bea
\sup_{n \geq 1, x \in \R^d} E[ H_{t,n}(x)^k] < \infty.
\lbl{060818}
\eea
By \eq{Htndef}, \eq{060818} and Fubini's theorem, 
\bea
E[\bar{H}_{t,n}] = |Q_n|^{-1}  \int_{Q_n} E H_{t,n}(x) dx.
\lbl{020808a}
\eea
By \eq{expest} and \eq{020808a}, using  the condition \eq{Qn3} on $Q_n$
it is straightforward to deduce that  $E\bar{H}_{t,n} \to EH_t(\0)$
as $n \to \infty$.

An alternative representation of $\bar{H}_{t,n}$ is as the sum
\bea
\bar{H}_{t,n} = |Q_n|^{-1} \sum_{(X,T) \in \Po_{t,n}} \xi((X,T),\Po_{t,n})
\lbl{060817b}
\eea 
where $\xi((X,T),\Po_{t,n})$ is integrated height 
of the vertically exposed parts 
(at time $t$) of the ball arriving at $(X,T)$. 
In other words, taking $J_{(X,T)}(x)$ to be the 
indicator function of the event that 
 the closure of the ball deposited above  location $X$ 
at time $T$ includes a point with  with coordinates
 $(x,H_{t,n}(x))$, we set
$$
\xi((X,T),\Po_{t,n}) = \int_{\R^d} H_{t,n}(x) J_{(X,T)}(x) dx.
$$
Using \eq{060817b}, we can obtain a law of
large numbers with almost sure convergence and convergence of 
the mean,
for $\bar{H}_{t,n}(x)$ by use of Theorem 3.2 of \cite{PY2}
(this is where we use conditions \eq{Qn4} and \eq{Qn5} on $Q_n$).
As mentioned above, the mean converges to
$E H_t(\0)$, so the limiting constant must be
$E H_t(\0)$. This demonstrates \eq{ctslln1}

  We now prove \eq{ctslln2}.
By expanding out the square in \eq{Vtn0}, we obtain the identity  
\bea 
V_{t,n}^2 =
 \left( |Q_n|^{-1} \int_{Q_n} H_{t,n}^2(x) dx \right) - \bar{H}_{t,n}^2,  
\lbl{Vtn}
\eea
so that by Fubini's theorem 
\bean
E[ V^2_{t,n}] 
+ [ E\bar{H}_{t,n}]^2
=  |Q_n|^{-1} \int_{Q_n} E[ H^2_{t,n}(x)]. 
\eean
Using  \eq{expest}, \eq{060818} and \eq{Qn3}, it is then straightforward to 
deduce that
\bea
E[ V^2_{t,n}] 
+ [ E\bar{H}_{t,n}]^2
\to 
E[(H_t(\0))^2] {\rm ~~~ as ~~} n \to \infty.
\lbl{060818a}
\eea

By \eq{Vtn}, we have
\bea
V_{t,n}^2 + \bar{H}_{t,n}^2 = |Q_n|^{-1} \sum_{(X,T) \in \Po_{t,n}}
\xi^{(2)}((X,T),\Po_{t,n})
\lbl{060815c}
\eea
where we set
$$
\xi^{(2)}((X,T),\Po_{t,n}) = \int_{\R^d} H_{t,n}(x)^2 J_{(X,T)}(x) dx.
$$
Hence we may
  use  Theorem 3.2 of \cite{PY2}
to deduce that the expression \eq{060815c} converges
to a deterministic limit,
 almost surely and with corresponding convergence of mean
(the `bounded $p$th moments condition' in \cite{PY2}
can be shown for this case by arguments in \cite{PY2}.)
Hence by \eq{060818a}, the almost sure 
limit of the expression \eq{060815c} is equal to 
$E[(H_t(\0))^2]$.
Combining this with \eq{ctslln1}, we obtain \eq{ctslln2}.  $\qed$ \\

For $t \geq 0$ and $u \geq 0$,
let $\tD_t$ denote the depth (i.e.  maximum height) of the dual
 continuum BD model at time $t$,
 and let $\tT(u) $ be the first passage time to depth $u$ of
 the dual continuum BD model.
 i.e. let
$$
 \tD_t :=  \sup_{x \in \R^d} ( \hat{H}_t(x)); ~~~ 
\tT(u) := \inf \{ t \geq 0: D_t \geq u\}.
$$

\begin{prop}
\lbl{thctsdual}
The distribution of $H_t(\0)$ is the same as that
of $\tD_t$.
\end{prop}

{\em Proof.}
Fix $t>0$, and consider the time-reversed space-time
Poisson process with time-coordinates transformed
by the mapping $(X,T) \mapsto (X,t-T)$.
Under this mapping, 
 any  path $((X_1,T_1), \ldots, (X_k,T_k))$
corresponds (by reversing the order of points)
to a path in the transformed Poisson process, namely
  $((X_k,t-T_k), \ldots, (X_1,t-T_1))$, the so-called {\em time-reversed
path}.
If the original path ends near $\0$,
the corresponding time-reversed path 
starts near 
 $\0$, and the height at $\0$  of the original path equals the
dual height over $X_1$ of the time-reversed path,
which is the maximal dual height of the time-reversed path.

By Lemma \ref{lempath} (ii), 
$H_t(\0) \geq u$ if and only if
there is a path in $\Po_t$ which ends near $\0$
with height at $\0$ of at least $u$ by time $t$, in which case
the corresponding time-reversed 
path has maximal height   at least $u$.
Hence, $H_t(\0)$ 
 is the maximal dual height 
of time-reversed 
paths in $\Po_t$  starting near $\0$. 
Hence by Lemma \ref{lempath} (iv), $H_t(\0)$ is the maximal depth
 at time $t$
for the continuum ballistic deposition
process driven by the transformed Poisson process,
using  initial profile $\hat{H}_0$.
Since the distributions of the original and transformed
Poisson processes are identical, 
 $H_t(\0)$ therefore
 has the same distribution as maximal depth 
 in the BD process generated by the original Poisson
process with initial profile $\hat{H}_0$. In other
words, it has the same distribution as $\tilde{D}_t$.
$\qed$ \\

The next result is a continuum analogue to Lemma 
\ref{lemmaxlim}.
\begin{lemm}
\lbl{lemmaxlim2}
There is a constant $\rho_3 \in (0,\infty)$ such that
\bea
\lim_{u \to \infty}
\frac{\tT(u)}{u} = \rho_3, ~~~a.s.
\label{060705}
\eea
and
\bea
\lim_{t \to \infty}
\frac{\tD_t}{t} = \rho_3^{-1}, ~~~a.s.
\label{060705b}
\eea
\end{lemm}
{\em Proof.}
First we show that $\tT(1)$ has finite second moment. 
Choose $\eps_1 >0$ such that
$F(\eps_1) \leq 1/2$, and $\eps_2 >0$ such that a ball of
radius $\eps_1$  in $\R^{d+1}$ contains a
rectilinear cube of side $\eps_2$,
with the same centre. Let $N_1(t)$ be the number of
Poisson arrivals in $\Po_t$ having  spatial coordinate 
with $\ell_\infty$ norm at most $\eps_2/2$ and having mark at least $\eps_1$.
Then 
$N_1(t)$ is Poisson with parameter at least $\eps_2 t/2$, 
and $\hat{H}_t(\0) \geq  \eps_2 N_1(t)$, so that
$$
P[ \tT(1) > t] 
\leq P[ \hat{H}_t(\0) < 1 ]
\leq P[ {\rm Po}(\eps_2 t/2) < 1/ \eps_2] =
 P\left[ \sum_{1 \leq i
\leq \lceil 1/\eps_2\rceil} e'_i >  t \right]
$$  
where $e'_1,e'_2,e'_3,\ldots$ are independent exponential 
random variables with mean $2/\eps_2$, representing
the inter-arrival times of a Poisson process
of rate $\eps_2/2$, and $\lceil x \rceil$ denotes the
smallest integer not less than $x$.
 Hence,
$$
E[\tT(1)^2 ] \leq  E \left[ \left(
\sum_{1 \leq i \leq  \lceil 1/\eps_2 \rceil} e'_i  \right)^2 \right] 
< \infty.
$$
We  assert that  $\tT(u)$ is 
distributionally subconvolutive, i.e. for $u, v \geq 0$ we have
\bea
F_{\tT(u+v)} \geq F_{\tT(u)} * F_{\tT(v)}.
\lbl{subconv2}
\eea
To see this, let $X(u) \in \R^d$ be chosen 
(in an arbitrary way if there is more than one choice)
so that $\hat{H}_{\tT(u)}(X(u)) \geq u$ (by definition such $X(u)$ exists).
 Let $(\hat{H}^*_s(x), x \in \R^d)_{s \geq 0}$ be a version of the BD process
with  initial profile
$$
\hat{H}^*_0(x)= 
 \left\{ \begin{array}{lr}
0 & \mbox{ if } x =X(u) \\ 
 -\infty & \mbox{ otherwise}
\end{array}
\right.    
$$
and driven by the Poisson process 
$\tau_{-\tT(u)} (\Po) \cap (\R^d \times (0,\infty))$,
where $\tau_{-t}$ denotes the shift operator on $\R^d \times \R$
mapping each point $(x,u) \in \R^d \times \R$ to $(x,u-t)$. 

Let $\tT^{*}(v)$ be the first time the process 
$(\hat{H}^*_s)_{s \geq 0}$ achieves
a depth at least $v$, i.e. 
$$
\tT^{*}(v) = \inf\{ t \geq 0: \sup_{x \in \R^d} 
(\hat{H}^*_t(x)) \geq v \}.
$$
Then $\tT^{*}(v)$ has the same distribution as $\tT(v)$,
and is independent of $\tT(u)$. Also the depth at time
 $\tT(u) + \tT^{*}(v)$ is at least $u+v$, i.e.
 $\tD_{\tT(u) + \tT^{*}(v)}
\geq u+v$, so that $\tT(u)  + \tT^{*}(v) \geq \tT(u+v)$.
Combining these facts gives us (\ref{subconv2}).

Since the variables $(\tT(u), u \geq 0)$ satisfy
\eq{subconv2} and are also
 monotonically increasing in $u$, we can apply 
 the Kesten-Hammersley theorem
(\cite{SW}, page 20) to
obtain the desired conclusion (\ref{060705}), with $0 \leq \rho_3 < \infty$.

The arguments to show that \eq{060705b} holds and $\rho_3 >0$
are just the same as the corresponding arguments in the proof of
Lemma \ref{lemmaxlim}.
%
$\qed$ \\

{\em Proof of Theorem \ref{ctsLLN}.}
The proof is entirely analogous to that of Theorem \ref{th020806}, now
using Propositions  \ref{LLN1} and  \ref{thctsdual} along with Lemma
 \ref{lemmaxlim2}.
$\qed$ \\ 

%
%
%

{\em  Proof of Theorem \ref{thctslb}.}
First we prove \eq{ctslower0}.
Let ${\cal F}_t$ be the $\sigma$-field generated
by all arrivals up to time $t-1$. 
Then $H_{t-1}(\0)$ is ${\cal F}_t$-measurable.
It is clear that there is a constant $\eps >0$ such
that
$$
P[ H_t (\0) = H_{t-1}(\0) | {\cal F}_t] \geq \eps ~~~ {\rm a.s.}
$$
and
$$
P[ H_t (\0) \geq H_{t-1}(\0) +1 | {\cal F}_t] \geq \eps ~~~ {\rm a.s.}
$$
and combining these two estimates shows that the conditional
variance ${\rm Var}[H_t|{\cal F}_t]$ is bounded
away from zero. From this we may deduce  \eq{ctslower0}.

The proof of \eq{ctslower} is similar to that of \eq{fplower}.
Now take $I_t$ to be the Lebesgue measure of the interface,
i.e., of the set of sites $x \in \R^d$ with $\hat{H}_t(x ) > -\infty$,
and take  ${\cal F}$ to be the $\sigma$-algebra generated by the
sequence of locations and marks of accepted arrivals in the dual continuum BD
process.

Then it is again the case that conditional on 
 ${\cal F}$,  
 the distribution of $T(u)$ is the sum of $M(u)$
independent exponentials $e_1,\ldots, e_{M(u)}$;
 this is because, given ${\cal F}$,
the $e_j$ are independent exponentials. Indeed, given the
positions of the first $j$ accepted arrivals, the
distribution of $e_j$ is exponential with mean $Y_j^{-1}$;
extra information about the location of this arrival and subsequent
arrivals does not affect its distribution. 

Using this information, the proof of \eq{ctslower} 
follows that of \eq{fplower} closely, and we omit further details.  $\qed$

Department of Mathematical Sciences, University
of Bath, Bath BA2 7AY, England.

 {\texttt masmdp@maths.bath.ac.uk} 

\end{document}